\documentclass[submission]{FPSAC2026}

\usepackage{tikz-cd}
\newtheorem{thm}{Theorem}
\newtheorem{lem}{Lemma}
\newtheorem{cor}{Corollary}

\newtheorem{prop}{Proposition}

\newtheorem{defn}{Definition}

\newcommand{\Q}{\mathbb{Q}}
\newcommand{\Z}{\mathbb{Z}}
\newcommand{\C}{\mathbb{C}}
\newcommand{\M}{\mathcal{M}}
\newcommand{\bbS}{\mathbb{S}}

\newcommand{\Mbar}{\overline{\mathcal{M}}}
\newcommand{\cat}[1]{\mathsf{#1}}
\newcommand{\Aut}{\mathrm{Aut}}

\newcommand{\calW}{\mathcal{W}}

\newcommand{\WRing}{\Lambda^{[2]}}
\newcommand{\ch}{\mathrm{ch}}
\newcommand{\Ind}{\mathrm{Ind}}
\newcommand{\pfunc}{\mathbf{P}}

\newcommand{\Part}{\mathsf{Part}}
\newcommand{\PPart}{\overline{\mathsf{Part}}}

\newcommand{\Partt}{\mathsf{Part}^{[2]}}
\DeclareMathOperator{\tcirc}{\tilde{\circ}}

\renewcommand{\-}{\text{-}}
\renewcommand{\P}{\mathbb{P}}

\usepackage{lipsum}

\title[P\'olya theory for moduli spaces of curves]{P\'olya enumeration, wreath product symmetric functions, and moduli spaces of curves}

\author[S. Kannan and T. Song]{Siddarth Kannan \thanks{\href{mailto:spkannan@mit.edu}{spkannan@mit.edu}. Siddarth Kannan is supported by NSF DMS-2401850.}\addressmark{1} \and Terry Dekun Song \thanks{\href{mailto:ds2016@cam.ac.uk}{ds2016@cam.ac.uk}. Terry Dekun Song is supported by a Cambridge Trust international scholarship.}\addressmark{2}}

\address{\addressmark{1}Department of Mathematics, Massachusetts Institute of Technology, Cambridge, MA \\ \addressmark{2}Department of Mathematics, University of Cambridge, Cambridge, UK}

\received{\today}

\abstract{We develop a calculus for $S_n$-equivariant Euler characteristics of moduli spaces of stable curves and stable maps. Our approach involves an enrichment of P\'olya's cycle index polynomial of a graph to a certain algebra $\WRing$ of wreath product symmetric functions. Building on foundational work of Macdonald, we prove that $\WRing$ may be viewed as the Grothendieck ring of the category of polynomial functors which map symmetric sequences of vector spaces to vector spaces. This interpretation gives rise to an action of $\WRing$ on the ordinary ring of symmetric functions $\Lambda$, which is described concretely in terms of Adams operations and skewing by power sums. This action lets us deduce appealing formulas, involving only ordinary symmetric functions, for generating functions of $S_n$-equivariant Euler characteristics.}


\usepackage[backend=bibtex]{biblatex}
\begin{document}

\maketitle

\section{Introduction}

For integers $g, n \geq 0$ such that $2g - 2 + n > 0$, let $\M_{g,n}$ denote the moduli space of genus-$g$ non-singular complex algebraic curves (=Riemann surfaces) which are compact and connected, together with a choice of $n$ distinct ordered points on the curve. Then $\M_{g,n}$ has the structure of a non-compact algebraic variety over $\C$, and it admits a compactification $\M_{g,n} \hookrightarrow \Mbar_{g,n}$, where $\Mbar_{g,n}$ is the Deligne--Mumford--Knudsen moduli space of stable $n$-pointed curves of genus $g$. Both $\M_{g,n}$ and $\Mbar_{g,n}$ are fundamental objects in geometry and topology, and the $S_n$-representations afforded by their rational cohomology groups have been the subject of intense study. Our first main theorem is a new formula which passes between the  generating functions \[ \cat{a} := \sum_{2g - 2 + n >0} \chi^{S_n} (\M_{g,n}) \cdot t^{g- 1} \quad \mbox{and}\quad \overline{\cat{a}} := \sum_{2g - 2 + n > 0} \chi^{S_n}(\Mbar_{g,n}) \cdot t^{g - 1} \]
for $S_n$-equivariant topological Euler characteristics of these moduli spaces. To be precise, the $S_n$-equivariant Euler characteristic of an $S_n$-variety $X$ is valued in the ring of symmetric functions $\Lambda$, and is defined by
\[ \chi^{S_n}(X) := \sum_i (-1)^i \mathrm{ch}_n(H^i(X;\Q)) \in \Lambda, \]
where $\ch_n(V) \in \Lambda$ denotes the Frobenius characteristic of a finite-dimensional $S_n$-representation.

Let us set some notation for our formula. We work over $\Q$, identifying
 \[ \Lambda = \Q[p_1, p_2, \ldots] \]
 where $p_i$ is the $i$th power sum symmetric function. We let $\Part$ denote the set of integer partitions. If $\mu \in \Part$, we write $\mu = (1^{\mu_1}, 2^{\mu_2}, \ldots )$, where $\mu_i$ is the number of parts of $\mu$ of size $i$. We set $p_\mu := \prod_{i} p_i^{\mu_i}$, and for $f \in \Lambda$ we define $\psi_\mu(f) = p_\mu \circ f$, where $\circ$ denotes plethysm of symmetric functions; each $\psi_\mu$ is a monomial in the so-called Adams operations on $\Lambda$.  We also write $p_\mu^\perp f$ for the operation of skewing with respect to $p_\mu$, given concretely by
\[ p_\mu^\perp f= \left({\prod_i i^{\mu_i}} \right)\frac{\partial^{\mu_1 + \mu_2 + \cdots}f}{\partial p_1^{\mu_1} \partial p_2^{\mu_2}\cdots}.\]
Our formula expands $\overline{\cat{a}}$ as an infinite sum, where each term is a monomial in expressions of the form $\psi_\mu(p_\lambda^\perp \cat{a})$ for $\mu, \lambda \in \Part$. Convenient indexing for the set of such monomials is provided by the set $\Partt$ of \textit{$2$-partitions}, defined by
\[ \Partt := \{ \Theta: \Part \to \Part \mid \Theta(\mu) = \varnothing \text{ for }|\mu|\gg 0 \}. \]
In other words, $\Partt$ is the set of finitely-supported functions $\Part \to \Part$.
Given $\Theta \in \Partt$, we define an operator $D_\Theta : \Lambda \to \Lambda$ by
\[D_\Theta(f) := \prod_{\mu \in \Part} \psi_{\Theta(\mu)}(p_\mu^\perp f).\]
Thus the symmetric functions $D_\Theta(f)$ for $\Theta \in \Partt$ index all monomials in terms of the form $\psi_\mu(p_\lambda^\perp f)$. We require the notation
\[||\Theta|| := \sum_{\mu \in \Part} |\mu| \cdot |\Theta(\mu)| \in \Z \]
for $\Theta \in \Partt$.

\begin{thm}\label{thm:mainthm}
    We have
    \[ \overline{\cat{a}} = \sum_{\Theta \in \Partt} O(\Theta) \cdot D_\Theta(\cat{a}) \cdot t^{||\Theta||/2}, \]
    where each coefficient $O(\Theta) \in \Q$ is the solution to a graph enumeration problem defined by (\ref{eqn:Otheta}) in \S \ref{subsec:graph_coeffs} below. 
\end{thm}
Note that $O(\Theta) = 0$ if $||\Theta||$ is odd, so $||\Theta||/2$ is an integer. Theorem \ref{thm:mainthm} also holds if we view the Euler characteristic as a virtual mixed Hodge structure; see \cite[\S1]{GraphEnumeration}.

\subsection{The coefficients $O(\Theta)$}\label{subsec:graph_coeffs}
Let $G$ be a finite graph, possibly containing loops or parallel edges. For each integer $i \geq 0$, set $\nu_i(G)$ to be the number of vertices of $G$ which have valence $i$. The sequence \[\nu(G) = (0^{\nu_0(G)}, 1^{\nu_1(G)}, \ldots)\]
which arises is an example of a \textit{generalized partition}, since there are possibly parts of size $0$. Let $\PPart$ denote the set of all generalized partitions. For each $\nu \in \PPart$, define
\begin{equation}\label{eqn:Snu}
\bbS_{\nu} := \prod_{i \geq 0} S_i \wr S_{\nu_i},
\end{equation}
where for integers $m, n \geq 0$ we have set $S_m \wr S_n := S_m^n \rtimes S_n$ for the wreath product of $S_m$ with $S_n$. For any graph $G$, there is a natural embedding
\begin{equation}\label{eqn:aut_graph_embedding}
\Aut(G) \hookrightarrow \bbS_{\nu(G)},
\end{equation}
well-defined up to conjugacy. The conjugacy classes of the groups $\bbS_\nu$ for $\nu \in \PPart$ have long been well-understood, essentially going back to Specht's dissertation; see \cite[Appendix I.B]{Macdonald}. The conjugacy classes are indexed by $2$-partitions.
\begin{lem}[Specht]\label{lem:Specht}
There is a natural bijection
\[ \bigcup_{\nu \in \PPart} \{\text{conjugacy classes of }\bbS_\nu \} \leftrightarrow \Partt. \]
\end{lem}

Using the correspondence of Lemma \ref{lem:Specht} and the embedding (\ref{eqn:aut_graph_embedding}), each $\Theta \in \Partt$ defines a (possibly empty) union of conjugacy classes $\Aut^{\Theta}(G) \subset \Aut(G)$, for any finite graph $G$. The rational numbers $O(\Theta)$ of Theorem \ref{thm:mainthm} are defined by
\begin{equation}\label{eqn:Otheta}
    O(\Theta) := \sum_{G \text{ connected graph}} \frac{|\Aut^\Theta(G)|}{|\Aut(G)|}.
\end{equation}
The above sum makes sense since $\Aut^\Theta(G) \neq \varnothing$ for only finitely many graphs $G$.

An alternative (and extremely elegant) formula to that of Theorem \ref{thm:mainthm} was previously found by Getzler--Kapranov \cite[Theorem 8.13]{GetzlerKapranov}. Our work is independent from theirs, and it is not clear how to derive either formula from the other. See \cite[\S 1.4.1]{GraphEnumeration} for an in-depth comparison of the two formulas. A beautiful expression for $\chi^{S_n}(\M_{g, n})$ has been given by Gorsky \cite{Gorsky}, so Theorem \ref{thm:mainthm} determines $\chi^{S_n}(\Mbar_{g,n})$ for all $g$ and $n$ in terms of graph sums. Similar graph sums, in the language of ``orbigraphs", appear via matrix integrals in work of Bini--Harer \cite{BiniHarer} on the non-equivariant Euler characteristics of these moduli spaces. Signed graph sums related to ours appear in \cite{CFGP} in the context of the weight-zero Euler characteristic of $\M_{g, n}$. A more thorough and conceptual explanation of related work is undertaken in our main article \cite{GraphEnumeration}.

For us, the most pleasing aspects of Theorem \ref{thm:mainthm} are the graph-theoretic interpretation of the coefficients $O(\Theta)$ and in the flexibility of the strategy developed to prove it. The strategy applies broadly to compute Euler characteristics of many related moduli spaces. Our second main theorem calculates the $S_n$-equivariant Euler characteristics of the moduli spaces of stable maps which are the foundation for Gromov--Witten theory.

\subsection{The moduli space of stable maps}\label{sec:stable_maps}The Kontsevich moduli spaces $\Mbar_{g,n}(\P^r, d)$ of degree-$d$ stable maps from $n$-pointed curves to the projective space $\P^r$ are of central interest in enumerative algebraic geometry, because their cohomology rings contain the information of the Gromov--Witten theory of $\P^r$. Prior to Theorem \ref{thm:stable_maps} below, calculations of $\chi^{S_n}(\Mbar_{g, n}(\P^r, d))$ were limited to the case $g = 0$ \cite{GetzlerPandharipande}, in which case the space is smooth and irreducible. When $g > 0$ these spaces are generally reducible, and have rather complicated singularities \cite{VakilsLaw}.



Our formula passes through torus localization \cite{graberpand, KontsevichTorusActions}, which leads us to consider the set $\Gamma_{\P^r, d}$ of isomorphism classes of tuples $(G, f, \delta)$ where
\begin{itemize}
\item $G$ is a finite connected graph;
\item $f: V(G) \to \{0, \ldots, r\}$ is a proper coloring of $G$ using at most $r + 1$ colors;
\item $\delta: E(G) \to \Z_{> 0}$ is a weighting by positive integers which satisfies
\[ \sum_{e \in E(G)} \delta(e) = d. \]
\end{itemize}
We can think of $\Gamma_{\P^r, d}$ as a groupoid, where the isomorphisms are required to respect the decorations. The positivity condition on $\delta$ implies that $\Gamma_{\P^r, d}$ is finite.

If $d = 0$, then $\Mbar_{g, n}(\P^r, 0) \cong \Mbar_{g, n} \times \P^r$, and if $r = 0$, then the target is a point and we recover $\Mbar_{g, n}$. Therefore we assume $r, d > 0$ and define
\[ \overline{\cat{a}}_{\P^r, d} := \sum_{g,n \geq 0} \chi^{S_n}(\Mbar_{g, n}(\P^r, d)) \cdot t^{g - 1}. \]
We then set
\[ \overline{\cat{a}}_{\dagger} := (h_1 + h_2)\cdot t^{-1} + \overline{\cat{a}}, \]
where $h_n \in \Lambda$ denotes the $n$th homogeneous symmetric function. The terms $h_1$ and $h_2$ can be interpreted as formally setting $\Mbar_{0, 1} = \Mbar_{0, 2} = \{\mathrm{pt}\}$. The embedding (\ref{eqn:aut_graph_embedding}) generalizes in a straightforward way to elements of $\Gamma_{\P^r, d}$.
\begin{thm}\label{thm:stable_maps}
We have
\[ \overline{\cat{a}}_{\P^r, d} = \sum_{\Theta \in \Partt} O_{\P^r, d}(\Theta) \cdot D_{\Theta}(\overline{\cat{a}}_{\dagger}) \cdot t^{||\Theta||/2}, \]
where
\[ O_{\P^r, d}(\Theta) = \sum_{(G, f, \delta) \in \Gamma_{\P^r, d}} \frac{|\Aut^{\Theta}(G)|}{|\Aut(G)|}. \]
\end{thm}

Both sums in Theorem \ref{thm:stable_maps} are finite, since $\Gamma_{\P^r, d}$ is finite. The finiteness is in stark contrast with Theorem \ref{thm:mainthm}. Theorem \ref{thm:stable_maps} makes many new calculations possible. See Table \ref{table:stable_maps_no_markings} for some sample calculations of $\Mbar_{g, 0}(\P^r, 3)$. We find it pleasing that while the variety $\Mbar_{g, n}(\P^r, d)$ governs the enumeration of algebraic curves in $\P^r$, its topological Euler characteristic is governed by the enumeration of $(r + 1)$-colored graphs. This perspective leads to the following structural result, which follows from properties of the chromatic polynomial.
\begin{cor}\label{cor: polynomiality}
There exists a polynomial $K_{g, n, d}(r) \in \Lambda[r]$ of degree $d$ such that
\[ \chi^{S_n}(\Mbar_{g, n}(\P^r, d)) = (r + 1) \cdot K_{g, n, d}(r + 1) \]
for all $r$.
\end{cor}
Corollary \ref{cor: polynomiality} suggests interesting stabilization properties for the cohomology of these moduli spaces as $r \to \infty$; in particular $\chi^{S_n}(\Mbar_{g, n}(\P^r, d))$ is determined, for all $r$, by the first $d$ values of $r$.

\begin{table}[h]
\centering
\def\arraystretch{1.5}
\begin{tabular}{|c|c|}
\hline
$g$ & $\chi \left(\Mbar_{g, 0}(\P^r, 3)\right)$                                                           \\ \hline
$0$ & 

$16\binom{r+1}{4}+21\binom{r+1}{3}+6\binom{r+1}{2}$                                 \\ \hline
$1$ & 
$216\binom{r+1}{4}+247\binom{r+1}{3}+55\binom{r+1}{2}$                                    \\ \hline
$2$ & $3160\binom{r+1}{4}+3342\binom{r+1}{3}+645\binom{r+1}{2}$ \\ \hline

$3$ & $44800\binom{r+1}{4}+45114\binom{r+1}{3}+8088\binom{r+1}{2}$ \\ \hline

$4$ & $630352\binom{r+1}{4}+613213\binom{r+1}{3}+104208\binom{r+1}{2}
$ \\ \hline
\end{tabular}

\caption{The topological Euler characteristic of $\Mbar_{g, 0}(\P^r, 3)$, for $g \leq 4$.}
\label{table:stable_maps_no_markings}
\end{table} 

Though the formulas in Theorems \ref{thm:mainthm} and \ref{thm:stable_maps} involve only Adams operations and skewing operators acting on ordinary symmetric functions, our proof passes through a certain algebra $\WRing$ of wreath product symmetric functions. We now turn to our combinatorial and representation-theoretic study of this ring, which should be of independent interest.

\section{Wreath symmetric functions and Schur--Weyl duality}
We will begin by recalling Schur--Weyl duality in the context of ordinary symmetric functions. Write $\cat{Vect}_{\Q}$ for the category of finite-dimensional vector spaces over $\Q$. 

\begin{defn}
    An \textbf{$\bbS$-module} or \textbf{symmetric sequence} over $\Q$ is a vector space $V \in \cat{Vect}_\Q$ together with a direct sum decomposition
    \[ V = \bigoplus_{n \geq 0} V_n, \]
    where each $V_n$ is an $S_n$-representation in $\cat{Vect}_\Q$. A morphism of $\bbS$-modules $V \to W$ is a sequence of $S_n$-equivariant linear maps $V_n \to W_n$ for all $n \geq 0$. We write $\bbS\-\cat{Vect}_\Q$ for the category of $\bbS$-modules over $\Q$. 
\end{defn}

Since an $\bbS$-module $V$ is finite-dimensional, we must have $V_n = 0$ for $n \gg 0$. The tensor product of two $\bbS$-modules is defined by
\begin{equation}\label{eqn:Stensorprod}
V \boxtimes W := \bigoplus_{n \geq 0} \left(\bigoplus_{k = 0}^{n} \Ind_{S_k \times S_{n - k}}^{S_n} V_k \otimes W_{n-k} \right).
\end{equation}
If we let $K_0(\bbS\-\cat{Vect}_\Q)$ denote the Grothendieck ring of $\bbS$-modules, we have
\begin{equation}\label{eqn:K_0_sym_seq}
    K_0(\bbS\-\cat{Vect}_\Q) \otimes_\Z \Q \cong \Lambda,
\end{equation}
by sending an $\bbS$-module $V$ to the sum $\sum_{n \geq 0} \ch_n(V_n)$. Schur--Weyl duality also identifies $\Lambda$ with the Grothendieck ring of \textit{polynomial functors} $\cat{Vect}_\Q \to \cat{Vect}_\Q$.
\begin{defn}\label{defn:pfun}
A functor $F: \cat{Vect}_\Q \to \cat{Vect}_\Q$ is called \textbf{polynomial} if for all $U, V \in \cat{Vect}_\Q$, the map
\[ \mathrm{Hom}_{\Q}(U, V) \to \mathrm{Hom}_{\Q}(F(U), F(V)) \]
induced by $F$ is a polynomial mapping of vector spaces. We write $\pfunc[\cat{Vect}_\Q, \cat{Vect}_\Q]$ for the category of polynomial functors $\cat{Vect}_\Q \to \cat{Vect}_\Q$, where the morphisms are natural transformations of functors.
\end{defn}

The following theorem is one version of the classical Schur--Weyl duality.

\begin{thm}[Schur--Weyl duality for $\bbS$-modules]\label{thm:SchurWeylClassical}
There is an equivalence of categories \[\bbS\-\cat{Vect}_\Q \to \pfunc[\cat{Vect}_\Q, \cat{Vect}_\Q],\] which takes an $\bbS$-module $V$ to the polynomial functor $F_V$ defined on a vector space $W$ by
\[ F_V(W) = \bigoplus_{n \geq 0} \mathrm{Hom}_{S_n}(V_n, W^{\otimes n}). \]
\end{thm}

Theorem \ref{thm:SchurWeylClassical} identifies the ring $\Lambda$ of symmetric functions with the Grothendieck ring of the category $\pfunc[\cat{Vect}_\Q, \cat{Vect}_\Q]$. The composition of polynomial functors descends to plethysm of symmetric functions.   

\subsection{$\bbS^{[2]}$-modules and polynomial functors of $\bbS$-modules}
Recall that to each $\nu \in \PPart$ we have associated in (\ref{eqn:Snu}) a finite group $\bbS_\nu$.
\begin{defn}
An \textbf{$\bbS^{[2]}$-module} is a vector space $V \in \cat{Vect}_\Q$ together with a direct sum decomposition
\[V = \bigoplus_{\nu \in \PPart} V_\nu, \]
where each $V_\nu$ is a representation of $\bbS_\nu$. A morphism of $\bbS^{[2]}$-modules is a collection of equivariant linear maps $V_\nu \to W_\nu$ for all $\nu$. We write $\bbS^{[2]}\-\cat{Vect}_\Q$ for the category of $\bbS^{[2]}$-modules.
\end{defn}


Define the algebra
\[ \WRing := \Q[p_1(\mu), p_2(\mu), \ldots \mid \mu \in \Part]. \]
We call $\WRing$ the ring of \textit{$2$-symmetric functions}. Each $p_i(\mu)$ should be thought of as a power sum indexed by the partition $\mu$. With some details omitted, Specht's classification of the conjugacy classes of $\bbS_\nu$ (Lemma \ref{lem:Specht}) gives an identification
\begin{equation}\label{eqn:grothendieck_S2}
K_0(\bbS^{[2]}\-\cat{Vect}_\Q) \otimes_\Z \Q \cong \WRing,
\end{equation}
where $K_0(\bbS^{[2]}\-\cat{Vect}_\Q)$ is the Grothendieck ring of $\bbS^{[2]}$-modules. The tensor product on $\bbS^{[2]}$-modules is defined similarly to (\ref{eqn:Stensorprod}). 

Note that Definition \ref{defn:pfun} admits a straightforward generalization to polynomial functors $\cat{A} \to \cat{Vect}_\Q$ whenever $\cat{A}$ is a subcategory of $\cat{Vect}_\Q$. Let $\pfunc[\bbS\-\cat{Vect}_\Q, \cat{Vect}_\Q]$ denote the category of polynomial functors $\bbS\-\cat{Vect}_\Q \to \cat{Vect}_\Q$. Building on foundational work of Macdonald \cite{MacDonald1980}, we prove the following theorem.
\begin{thm}[Schur--Weyl duality for $\bbS^{[2]}$-modules]\label{thm:upgradedSchurWeyl} There is an equivalence of categories
\[ \bbS^{[2]}\-\cat{Vect}_\Q \to \pfunc[\bbS\-\cat{Vect}_\Q, \cat{Vect}_\Q],  \]
which takes a $\bbS^{[2]}$-module $V$ to the functor $F_V$ determined by
\[ F_V\left(\bigoplus_{n \geq 0} W_n\right) = \bigoplus_{\nu \in \PPart} \mathrm{Hom}_{\bbS_\nu}\left(V_\nu, \bigotimes_{i \geq 0} W_i^{\otimes \nu_i} \right). \]
\end{thm}

We will now explain how Theorem \ref{thm:upgradedSchurWeyl} gives rise to an action
\[\WRing \times \Lambda \to \Lambda \]
which is the main tool for deriving the formulas of Theorems \ref{thm:mainthm} and \ref{thm:stable_maps}. 
\subsection{The action of $\WRing$ on $\Lambda$}
The composition of functors
\[ \pfunc[\bbS\-\cat{Vect}_\Q, \cat{Vect_\Q}] \times \pfunc[\cat{Vect}_\Q, \bbS\-\cat{Vect_\Q}] \to \pfunc[\cat{Vect}_\Q, \cat{Vect}_\Q] \]
descends via taking $K_0(\-) \otimes_{\Z} \Q$ to an action
\[ \WRing \times K_0(\pfunc[\cat{Vect}_\Q, \bbS\-\cat{Vect}_\Q])\otimes_\Z \Q\to \Lambda. \]
It is not too difficult, using the classical Schur--Weyl duality (Theorem \ref{thm:SchurWeylClassical}), to show that in fact 
\[ K_0(\pfunc[\cat{Vect}_\Q, \bbS\-\cat{Vect}_\Q])\otimes_\Z \Q \cong \Lambda \otimes_\Q \Lambda.\] Therefore composition of functors gives rise to a map
\begin{equation}\label{eqn:total_action}
    \WRing \times \Lambda \otimes_\Q \Lambda\to \Lambda.
\end{equation}
Now recall that the ring $\Lambda$ is a Hopf algebra, with coproduct $\Delta: \Lambda \to \Lambda \otimes_\Q \Lambda$ determined by $p_i \mapsto 1 \otimes p_i + p_i \otimes 1$ for all $i$. By applying (\ref{eqn:total_action}) to this coproduct, we obtain an action
 \[ \WRing \times \Lambda \to \Lambda, \]
 for which we use the notation
 \[(g, f) \mapsto g \tcirc f. \]
The action $\tcirc$ is described in terms of Adams operations and skewing in $\Lambda$.
 \begin{lem}\label{lem:tcirc}
The action $\tcirc$ is uniquely characterized by the following properties:
\begin{enumerate}
\item for all $f \in \Lambda$, we have $p_n(\mu) \tcirc f = \psi_n(p_\mu^\perp f)$, and
\item for fixed $f \in \Lambda$ the map $\WRing \to \Lambda$ by $g \mapsto g \tcirc f$ is an algebra homomorphism.
\end{enumerate}
\end{lem}
We now interpret $\tcirc$ representation-theoretically.
\subsection{Representation-theoretic interpretation of the action}
In view of (\ref{eqn:grothendieck_S2}), there is a direct sum decomposition
\[ \WRing = \bigoplus_{\nu\in\PPart}\mathrm{Cl}_\Q(\bbS_\nu),  \]
where $\mathrm{Cl}_\Q(\bbS_\nu)$ is the space of $\Q$-valued class functions on $\bbS_\nu$; this is the usual identification between class functions and virtual representations of a group. Via Lemma \ref{lem:Specht}, the $\Q$-algebra $\WRing$ can be naturally identified with the free vector space spanned by $2$-partitions. Given a $2$-partition $\Theta$, define
\[ p_\Theta := \prod_{\mu \in \Part} p_{\Theta(\mu)}(\mu) \in \WRing, \]
where we have set $p_{\lambda}(\mu):= \prod_{i > 0} p_{i}(\mu)^{\lambda_i}$ for $\lambda \in \Part$.

\begin{defn}
    For a finite-dimensional $\bbS_\nu$-representation $V_\nu$, the \textbf{Frobenius characteristic} $\ch_{\bbS_\nu}(V_\nu) \in \WRing$ is defined by
    \[\ch_{\bbS_\nu}(V_\nu) = \frac{1}{|\bbS_\nu|}\sum_{\boldsymbol{\tau} \in \bbS_\nu} \mathrm{Tr}(\boldsymbol{\tau}|V_{\nu})\cdot p_{\Theta^{\boldsymbol{\tau}}}, \]
    where $\Theta^{\boldsymbol{\tau}}$ is the $2$-partition encoding the cycle-type of $\boldsymbol{\tau}$ as in Lemma \ref{lem:Specht}. If $V$ is an $\bbS^{[2]}$-module, we define
    \[ \ch_{\bbS^{[2]}}(V) := \sum_{\nu \in \PPart} \ch_{\bbS_\nu}(V_\nu). \]
\end{defn}
Suppose now that we are given an $\bbS$-module $W$. From $W$ and an integer $i \geq 0$ we can make a new $\bbS$-module $W(i)$ by
\[ W(i)_n = \mathrm{Res}^{S_{i + n}}_{S_i \times S_n} W_{i + n}. \]
Each $W(i)$ is naturally an $S_i$-object in the category $\bbS\-\cat{Vect}_\Q$. This makes the $j$th tensor power $W(i)^{\boxtimes j}$ into a $S_i \wr S_j$-representation in $\bbS\-\cat{Vect}$. For $\nu \in \PPart$, we define the $\bbS$-module
\[ W^{\otimes\nu}:= \bigoplus_{i \geq 0} W(i)^{\boxtimes \nu_i}, \]
which is an $\bbS_\nu$-representation in $\bbS\-\cat{Vect}$. Now if \[V = \bigoplus_{\nu \in \PPart} V_\nu \]
is an $\bbS^{[2]}$-module, we get an endofunctor
\[ T_{V} : \bbS\-\cat{Vect} \to \bbS\-\cat{Vect} \]
by
\[ T_{V}(\calW) = \bigoplus_{\nu \in \PPart} \mathrm{Hom}_{\bbS_\nu}(V_\nu, W^{\otimes\nu}). \]
The functor $T_V$ gives a categorification of the action $\WRing \times \Lambda \to \Lambda$.
\begin{prop}\label{prop:PPaction}
For an $\bbS$-module $W$ and an $\bbS^{[2]}$-module $V$, we have
\[ \ch_{\bbS}(T_{V}(W)) = \ch_{\bbS^{[2]}}(V) \tcirc  \ch_{\bbS}(W). \]
\end{prop}
To apply our theory to the moduli spaces of curves, we use Proposition \ref{prop:PPaction} in the context of graph automorphism groups.

\section{The P\'olya--Petersen character of a graph}

\begin{figure}[h]
    \centering
    \includegraphics[scale=.75]{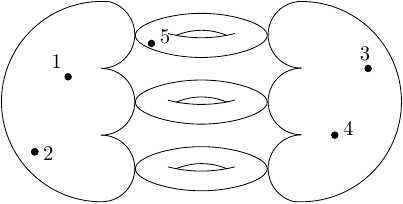}
    \caption{A point of $\Mbar_{5, 5}$ with dual graph given in Figure \ref{fig:graph_exmp}}
    \label{fig:stable_curve_exmp}
\end{figure}
\begin{figure}
        \centering
        \includegraphics[scale=1]{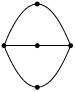}
        \caption{The dual graph of the nodal algebraic curve in Figure \ref{fig:stable_curve_exmp}}    
        \label{fig:graph_exmp}
    \end{figure}

A basic geometric property of the moduli space $\Mbar_{g,n}$ is that it admits a stratification by the \textit{dual graph} of the underlying curve. The dual graph of a nodal algebraic curve $C$ has a vertex for each irreducible component of $C$, and an edge between two vertices for each node shared by the corresponding components (this includes the possibility of loops). Note that this is a coarser notion of a dual graph than is usual for stable algebraic curves. See Figure \ref{fig:stable_curve_exmp} for an example of a point of the moduli space $\Mbar_{5,5}$ which has the dual graph given in Figure \ref{fig:graph_exmp}.

Now if $G$ is any finite graph, recall that there is a natural embedding (\ref{eqn:aut_graph_embedding})
\[ \Aut(G) \hookrightarrow \bbS_{\nu(G)}, \]
well-defined up to conjugacy.
\begin{defn}
For a finite graph $G$, the \textbf{P\'olya--Petersen character} of $G$ is the element
\[ \zeta_G := \mathrm{ch}_{\bbS_{\nu(G)}}\left( \Ind_{\Aut(G)}^{\bbS_{\nu(G)}} \mathrm{triv} \right) \in \WRing. \]
\end{defn}

For example, consider the graph $G$ given in Figure \ref{fig:graph_exmp}. In this case $\nu(G) = (2^3,3^2)$, so
\[\bbS_{\nu(G)} = (S_2 \wr S_3) \times (S_3 \wr S_2), \]
and $\Aut(G)$ is an extension of $S_3$ by $\Z/2\Z$. We may view the $\Z/2\Z$ factor as coming from a horizontal flip and the elements of $S_3$ as permutations of the bivalent vertices and their adjacent edges.

We may compute the P\'olya--Petersen character of $G$ as
\begin{align*}
\zeta_G = \frac{1}{12}\left( \underbrace{p_1(1^2)^3 p_1(1^3)^2}_{\text{identity}} \right. &+ \underbrace{3p_1(1^2) p_1(1^12^1)p_2(1^2)}_{2\text{-cycles without flip}} + \underbrace{2p_1(3)^2p_3(1^2)}_{3\text{-cycles without flip}} \\& + \left.\underbrace{p_1(2^1)^3 p_2(1^3)}_{\text{flip}} + \underbrace{3 p_1(2^1) p_2(1^2) p_2(1^3)}_{2\text{-cycles with flip}} + \underbrace{2 p_2(3^1) p_3(2^1) }_{3\text{-cycles with flip}} \right).
\end{align*}

The key application of the P\'olya--Petersen character to Theorem \ref{thm:mainthm} is the following Proposition, which follows from the representation-theoretic interpretation of the action $\tcirc$ of $\WRing$ on $\Lambda$. For each finite graph $G$ and integers $g, n$ satisfying $2g - 2 + n > 0$, define
\[ \Mbar_{g, n}(G) \subset \Mbar_{g, n} \]
to be the subspace of curves with dual graph equal to $G$.

\begin{prop}\label{prop: key_prop}
    If we set
    \[ \cat{a}_G := \sum_{2g - 2 + n > 0} \chi^{S_n}(\Mbar_{g, n}(G)) \cdot t^{g - 1}, \]
    then
    \[ \cat{a}_G = (\zeta_G \tcirc \cat{a}) \cdot t^{|E(G)|}. \]
\end{prop}
Theorem \ref{thm:mainthm} is deduced from Proposition \ref{prop: key_prop} by summing both sides over all finite connected graphs and using the explicit description of $\tcirc$ in Lemma \ref{lem:tcirc}. The power $||\Theta||/2$ appearing in the formula is equal to $|E(G)|$ for any graph $G$ with $\Aut^{\Theta}(G) \neq \varnothing$.

We conclude by remarking on the choice of name ``P\'olya--Petersen character.'' In the context of his enumeration theorem, P\'olya defined the cycle index polynomial of a group $H$ acting on a finite set $X$ \cite{Polya}. In modern language, his polynomial is the symmetric function defined by
\[ z_{H, X}^{\textrm{P\'ol}} := \ch_{|X|}(\Ind_{H}^{S_{|X|}}\mathrm{triv}) \in \Lambda. \]
P\'olya's cycle index for a finite graph $G$ is recovered when $H = \Aut(G)$ and $X = V(G)$ is the set of vertices. The P\'olya--Petersen character $\zeta_G$ specializes to the cycle index $z_{\Aut(G), V(G)}^{\textrm{P\'ol}}$ under the $\Q$-algebra map $\WRing \to \Lambda$ defined by $p_n(\mu) \mapsto p_n$ for all $n$ and $\mu$.

Our approach owes an intellectual debt to Petersen \cite{semiclassicalremark}, who used similar techniques to calculate $\chi^{S_n}(\Mbar_{1,n})$. For genus one curves, dual graphs are not much more complicated than cycles, so Petersen was led to consider the character $\zeta_G$ in the special case where $G$ is a $k$-cycle. In this case $\bbS_{\nu(G)}$ is the $k$th hyperoctahedral group, so in his work $\zeta_G$ is a type-$B$ Frobenius characteristic in the ring $\Lambda(S_2)$ of $S_2$-wreath symmetric functions of Macdonald. We used similar ideas to compute $\chi^{S_n}(\Mbar_{1,n}(\P^r, d))$ quite explicitly by obtaining closed formulas for all the relevant graph sums \cite[Theorem A]{genusonechar}.








\acknowledgements{
We are grateful to Jonas Bergstr\"om, Sarah Brauner, Samir Canning, Melody Chan, Dan Petersen, Oscar Randal-Williams, and Dhruv Ranganathan for helpful conversations and correspondence. We also thank the reviewers for their careful reading of this abstract.
}

\printbibliography

\end{document}